\documentclass[leqno]{gtart}
\usepackage{amsmath}
\usepackage{amssymb}
\usepackage{amsthm}
\usepackage{amsfonts}
\usepackage{eucal}
\usepackage{xspace}
\usepackage{tabularx}
\usepackage{graphicx}
\input xy

\xyoption{all}

\usepackage{graphics}

\usepackage[all]{xy}

\NoCompileMatrices

\newcommand{\field}[1]{\ensuremath{\mathbb{#1}}}

\newcommand{\F}{\field{F}\xspace}

\newcommand{\N}{\field{N}\xspace}

\newcommand{\s}{\field{S}\xspace}

\newcommand{\Z}{\field{Z}\xspace}

\newcommand{\cuadro}{\hfill{$\sqcap \kern-8pt \sqcup $}}

\newcommand{\fy}{\varphi}
\newcommand{\en}{\longrightarrow}
\newcommand{\hook}{\hookrightarrow}
\newcommand{\onto}{\twoheadrightarrow}

\newcommand{\Mod}{{\mathcal M}{\rm od}}
\newcommand{\ktop}{\text{k-}{\mathcal T}{\rm op}}
\newcommand{\Gktop}{G\text{-k-}{\mathcal T}{\rm op}}
\newcommand{\Top}{{\mathcal T}{\rm op}}
\newcommand{\Set}{{\mathcal S}{\rm et}}
\newcommand{\Gset}{{G {\text{-}} {\mathcal S}{\rm et}}}
\newcommand{\Gtop}{{G {\text{-}} {\mathcal T}{\rm op}}}
\newcommand{\kmod}{{k {\text{-}} {\mathfrak Mod}}}
\newcommand{\Rmod}{{R {\text{-}} {\mathfrak Mod}}}

\newcommand{\EF}{{\mathcal F}}
\newcommand{\ha}{{\mathcal H}}
\newcommand{\ka}{{\mathcal K}}
\newcommand{\ga}{{\mathcal G}}
\newcommand{\Topab}{{\mathcal T}{\rm opab}}

\newcommand{\EFG}{{\mathcal F}^G}

\newcommand{\OEFG}{\overline{\mathcal F}^G}
\newcommand{\FG}{{F^G}}

\newcommand{\sing}{{\mathcal S}}

\DeclareMathOperator{\id}{id}

\DeclareMathOperator{\colim}{\mathop{\rm colim}}

\newtheorem{thm}{Theorem}[section]
\newtheorem{cor}[thm]{Corollary}
\newtheorem{lem}[thm]{Lemma}
\newtheorem{prop}[thm]{Proposition}

\newcommand{\qedsymb}{\mbox{
}~\hfill~{\rule{2mm}{2mm}}}
\newenvironment{pf}{\begin{trivlist}
\item[\hspace{\labelsep}{\it \noindent Proof: }]
}{\qedsymb\end{trivlist}}

\theoremstyle{definition}

\newtheorem{defs}[thm]{Definition}
\newtheorem{remark}[thm]{\sc Remark}

\newtheorem{examp}[thm]{\sc Example}
\newtheorem{examps}[thm]{\sc Examples}

\def\mapright#1{\smash{\mathop{\hbox to
40pt{\rightarrowfill}}\limits^{#1}}}

\def\supsetnoteq{\mathbin{\hbox{$\supseteq \joinrel \hskip-8pt
\lower3pt
                 \hbox{$\scriptscriptstyle /$}\ $}}}

\def\pieams#1{\begingroup \def \protect
    {\noexpand \protect \noexpand }\xdef \@thanks {\@thanks
    \protect \footnotetext [\the \c@footnote
    ]{\hbox{\hskip-19truept}#1}}\endgroup}

\begin{document}

\title{Equivariant Dold-Thom topological groups}

\author{Marcelo A. Aguilar \&  Carlos Prieto\footnote{Corresponding
  author, Phone: ++5255-56224489, Fax ++5255-56160348. The authors were partially supported by PAPIIT grant IN101909.}}

\address{Instituto de Matem\'aticas\\ Universidad Nacional Aut\'onoma
  de M\'exico\\ 04510 M\'exico, D.F.,
  Mexico}

\email{marcelo@math.unam.mx, cprieto@math.unam.mx}

\begin{abstract}
Let $M$ be a covariant coefficient system for a finite group $G$.  In this paper,
we analyze several topological abelian groups, some of them new, whose homotopy groups are isomorphic to the Bredon-Illman $G$-equivariant
homology theory with coefficients in $M$. We call these groups equivariant Dold-Thom topological groups and we show that they are unique up
to homotopy.  We use one of the new groups to prove that the Bredon-Illman homology satisfies the infinite-wedge axiom and to make some
calculations of the 0th equivariant homology
\end{abstract}




\maketitle

\pagenumbering{arabic}\setcounter{page}{1}

\setcounter{section}{-1}
\section{Introduction}\label{sec0}

Given a finite group $G$, a covariant coefficient system  $M$ for $G$,
and a pointed $G$-set $S$, we defined in \cite{mackey-3} an abelian
group $F^G(S,M)$. Using this construction, we associate to a pointed $G$-space $X$
a simplicial abelian group $\FG(\sing(X),M)$. Its
geometric realization, denoted by $\FG(X,M)$, is a topological group whose
homotopy groups are isomorphic to the reduced Bredon-Illman $G$-equiv\-ariant
homology of $X$ with coefficients in $M$.  Given a $G$-equivariant ordinary covering map
$p:E\en X$, we also defined a continuous transfer $t^G_p:\FG(X^+,M)\en \FG(E^+,M)$, that induces a transfer
in equivariant homology, when $M$ is a Mackey functor.

In \cite{bihomol} we showed that when the $G$-space $X$ is a strong $\rho$-space (e.g. a $G$-simplicial complex
or a finite dimensional countable locally finite $G$-CW-com\-plex),
there is a topology in the abelian groups $F^G(X^\delta,M)$, where $X^\delta$ stands
for the underlying set of $X$. With this topology, $\FG(X,M)$ is a topological group,
which we denote by $\F^G(X,M)$.  This is a smaller group than the group
$\FG(X,M)$ mentioned above. It has the property that its homotopy groups are
isomorphic to the reduced Bredon-Illman $G$-equivariant homology of $X$ with coefficients
in $M$, provided that $M$ is a homological Mackey functor. Furthermore, we proved in \cite{G-transrami} that these topological groups admit
a continuous transfer for $G$-equivariant ramified covering maps, whose total space and base space are strong $\rho$-spaces.

In this paper we present a different topology for the abelian group $F^G(X^\delta,M)$ and we denote
the resulting topological group by $\EFG(X,M)$.  We prove that for any pointed $G$-space $X$ of the
homotopy type of a $G$-CW-complex, and for any coefficient system, the homotopy groups of
$\EFG(X,M)$ are isomorphic to the reduced Bredon-Illman $G$-equivariant homology of $X$
with coefficients in $M$.  The assumptions on $X$ and $M$ are much weaker than those needed to define $\F^G(X,M)$.
However, in such a generality, it does not seem possible to construct a continuous transfer even for
ordinary covering $G$-maps.

These new topological groups $\EFG(X,M)$ can be used to prove the infinite wedge axiom
for the Bredon-Illman homology.  They also allow us to make some calculations of the 0th
homology groups of a $G$-space $X$.

Those topological groups whose homotopy groups are
isomorphic to the reduced Bredon-Illman $G$-equivariant homology of $X$ with coefficients
in $M$ will be called \emph{Dold-Thom topological groups}.
Hence $\FG(X,M)$, $\F^G(X,M)$, and $\EFG(X,M)$ are Dold-Thom topological groups.
Furthermore, all these groups are algebraically subgroups of other topological groups,
so they also have another natural topology, namely the subspace topology. These topological
groups will be denoted by
$\overline F^G(X,M)$, $\overline \F^G(X,M)$, and $\OEFG(X,M)$.
The first two were studied in \cite{bihomol}. In this paper we prove that these groups
are isomorphic to the former, if the coefficient system $M$ takes values on $k$-modules,
where $k$ is a field of characteristic $0$ or a prime $p$ that does not divide the order of $G$.
We also show that the Dold-Thom topological groups are unique up to homotopy. Other examples of Dold-Thom topological groups were constructed by Dold and Thom \cite{DT}, McCord \cite{mccord}, Lima-Filho \cite{Lima}, dos Santos \cite{santos}, and Nie \cite{nie}.

The paper is organized as follows.  In Section \ref{sec1} we give the basic definitions that are
needed.  Then in Section \ref{sec2} we define the new topological groups $\EFG(X,M)$ and
$\OEFG(X,M)$, and we prove that the former is a Dold-Thom topological group.  Then in Section
\ref{sec3} we compare the new topological groups with the previously defined ones. In
Section \ref{sec4} we analyze the case of coefficients in $k$-$\Mod$, with the field $k$
as explained above. In Section \ref{sec4.5} we prove the wedge axiom and we compute the 0th equivariant homology in some cases.
Finally in Section \ref{sec5}, we study the Dold-Thom topological groups
and we show that two Dold-Thom topological groups, which are locally connected and
have the homotopy type of a CW-complex, are homotopy equivalent.

\section[Preliminaries]{\large\sc Preliminaries}\label{sec1}

We shall work in the
category of k-spaces, which will be denoted by $\ktop$.  We understand by a k-{\em space} a
topological space $X$ with the property that a set $W\subset X$ is
closed if and only if $f^{-1}W\subset Z$ is closed for any
continuous map $f:Z\en X$, where $Z$ is any compact Hausdorff
space (see \cite{may-2,vogt}).  Given any space $X$, one can clearly associate to it a k-space
$k(X)$ using the condition above, which is weakly homotopy equivalent to $X$.
The product of two spaces $X$ and $Y$
in this category is $X\times Y = k(X\times_{\text{top}} Y)$, where
$X\times_{\text{top}} Y$ is the usual topological product.  Two important properties of this category
are that \emph{if $p:X\onto X'$ is an identification and $X$ is a {\rm k}-space, then $X'$ is a
{\rm k}-space}, and \emph{if $p:X\onto X'$ and $q:Y\onto Y'$ are identifications, then
$p\times q : X\times Y \onto X'\times Y'$ is an identification too}.  If $X$ is a k-space, we shall
say that $A\subset X$ has the \emph{relative {\rm k}-topology} (in $\ktop$) if $A = k(A_{\text{rel}})$, where
$A_{\text{rel}}$ denotes $A$ with the (usual) relative topology in $\Top$.  This topology is characterized
by the following property: \emph{Let $Y$ be a {\rm k}-space. Then a map $f:Y\en A$ is continuous
if and only if the composite $i\circ f:Y\en X$ is continuous, where $i:A\hook X$ is the inclusion} (see \cite{vogt}).

In what follows,
we shall denote by $\Gtop_*$ the category of topological pointed $G$-spaces
such that $G$ acts trivially on the base point, or correspondingly $\Gktop_*$.
$\Topab$ will denote the
category of topological abelian groups in the category of k-spaces.
Recall that a {\em  covariant coefficient system} is a covariant functor
   $M:{\mathcal O}(G)\en \Rmod$, where ${\mathcal O}(G)$ is the category of
$G$-orbits $G/H$, $H\subset G$, and $G$-functions $\alpha:G/H\en G/K$, and $R$ is
a commutative ring.   A
particular role will be played by the $G$-function $R_{g^{-1}}:G/H\en
  G/gHg^{-1}$,
given by {\em right translation} by $g^{-1}\in G$, namely
$$R_{g^{-1}}(aH)=aHg^{-1}=
ag^{-1}(gHg^{-1})\,.$$
We shall often denote $aH$ by $[a]_H$.  Observe
that if $X$ is a $G$-set and $x\in X$, then
the canonical bijection $G/G_x \en G/G_{gx}$ is precisely $R_{g^{-1}}$.
Here $G_x$ denotes the {\em isotropy subgroup} of $x$, namely, the
maximal subgroup of $G$ that leaves $x$ fixed.

Let $S$ be a pointed $G$-set (where the
base point $x_0$ remains fixed under the action of $G$) and $M$ a covariant coefficient system.
In \cite{mackey-3} we defined an abelian group $F(S,M)$ as follows. Consider the following union
$$\widehat M = \bigcup_{H\subset G} M(G/H)\,.$$
Then
$$F(S,M) = \{u:S\en \widehat M\mid  u(x)\in M(G/G_x)\,,\,\, u(x_0) = 0\,,\,\, $$
$$\text{and } u(x) = 0 \text{ for almost every }
x\in X\}\,.$$

Indeed, this group $F(S,M)$ is an $R$-module, whose structure is given by
$$(r\cdot u)(x) = ru(x)\in M(G/G_x)\,.$$
It has as canonical generators the functions
$$lx:S\en \widehat M\quad\text{ given by }\quad lx(x') = \begin{cases}
l & \text{ if } x' = x\,, \\
0 & \text{ if } x'\ne x\,,
\end{cases}$$
where $x\in S$, $x\ne x_0$, and $l\in M(G/G_x)$.  The group $F(S,M)$ is a functor
of $S$ as follows.  Let $f:S\en T$ be a pointed $G$-function.  Then we define
$f_*:F(S,M)\en F(T,M)$ on generators by
$$f_*(lx) = M_*(\widehat f_x)(l)f(x)\,,$$
that is, the homomorphism whose value on $y$ is $M_*(\widehat f_x)(l)$ if
$y = f(x)$ and $0$ otherwise, where $\widehat f_x:G/G_x\onto G/G_{f(x)}$ is
the canonical surjection.

There is an action of $G$ on $F(S,M)$ given on generators by
$$g\cdot (lx) = M_*(R_{g^{-1}})(l)(gx)\,.$$
Then one can consider the submodule $F(S,M)^G$ of
fixed points under the $G$-action.  Let $f:S\en T$ be as above.
Since $f_*$ is clearly a $G$-homomorphism,
it restricts to a homomorphism between the submodules of fixed points, which we
denote by
$$\overline f^G_*:\overline F^G(S,M)\en \overline F^G(T,M)\,.$$
This makes $\overline F^G(-,M)$ into a functor $\Gset_*\en \Rmod$

On the other hand, there is a surjective homomorphism $\beta_S:F(S,M)\onto F(S,M)^G$
given by
$$\beta_S(lx) = \gamma^G_x(l) = \sum_{[g]\in G/G_x} M_*(R_{g^{-1}})(l)(gx)\,,$$
that is, essentially taking the sum over the orbit of $x$.  One can now use this to
give a different functorial structure on $F(S,M)^G$, defining for a $G$-function
$f:S\en T$ and a generator $\gamma^G_x(l)$,
$$f^G_*(\gamma^G_x(l)) = \gamma^G_{f(x)}M_*(\widehat f_x)(l)\,.$$
This makes $F^G(-,M)$\label{functor} into a functor $\Gset_*\en \Rmod$.

These three functors are related by the commutativity of the following diagram:
\setcounter{thm}{1}
\begin{equation}\label{naturality}
\xymatrix{
S\ar[d]^f & & \overline F^G(S,M)\ar@{^(->}[r]^-{\iota_S}\ar[d]_-{\overline f^G_*} &
F(S,M)\ar@{->>}[r]^-{\beta_S}\ar[d]^-{f_*} & F^G(S,M)\ar[d]^-{f^G_*} \\
T & & \overline F^G(T,M)\ar@{^(->}[r]_-{\iota_T} &
F(S,M)\ar@{->>}[r]_-{\beta_T} & \,\,F^G(S,M)\,.}
\end{equation}

We shall use these groups to define different topological groups in the next section.

\section[Topological groups and coefficient systems]{\large\sc Topological groups and coefficient systems}\label{sec2}

First recall the following construction. Let $X$ be a topological space
and $L$ an $R$-module.
Then we have the $R$-module $F(X,L)=\{u:X\en L\mid u(x_0) = 0,
\text{ and } u(x) = 0 \text{ for almost every } x\in X \}$.
Following \cite{G-homol} we have that this $R$-module can be
topologized as follows (see also \cite{mccord}). There is a surjective function
$$\mu: \coprod_k (L\times X)^k\onto F(X,L)$$
given by $(l_1, x_1;\dots;l_k,x_k) \mapsto l_1 x_1+\cdots+l_kx_k$.
Then $F(X,L)$ has the identification topology.

Given a pointed $G$-space $X$, we denote by $X^H$ the pointed subspace of elements of
$X$ that remain fixed under the group elements of $H$.

\begin{defs}\label{garigolF} Let $M$ be a covariant coefficient system and let $X$ be any
  pointed $G$-space. For each subgroup $H\subset G$ consider the topological
 group $F(X^H,M(G/H))$ as defined above. Define $p_H:F(X^H,M(G/H))\en F(X,M)$
 by $p_H(lx) = M_*(q_{H,x})(l)$, where  $x\in X^H$, $l\in M(G/H)$, and $q_{H,x}:G/H\onto G/G_x$,
 is the canonical projection.
 Now take the homomorphism $$p_X: \prod_{H\subset G} F(X^H,M(G/H))\onto F(X,M)$$
 given by $p_X((l_Hx_H)_{H\subset G}) = \sum_{H\subset G} p_H(l_Hx_H)$,
where the product has the product topology of k-spaces.  Given any generator
$lx\in F(X,M)$, $lx$ can be seen also as an element in $F(X^{G_x},M(G/G_x))$.  Therefore
$p_X$ is surjective. Give the $R$-module $F(X,M)$ the identification topology
induced by $p_X$.  We obtain a k-space, which we denote
by $\EF(X,M)$.  By giving $\overline F^G(X,M)$ the relative
{\rm k}-topology we obtain another k-space, which we denote by $\OEFG(X,M)$.  Moreover, by
taking on $F^G(X,M)$ the identification topology given by the epimorphism $\beta_X$,
we obtain another k-space, which we denote by  $\EFG(X,M)$.
\end{defs}

\begin{prop}\label{garigoltop} The groups $\EF(X,M)$, $\OEFG(X,M)$, and $\EFG(X,M)$ are
topological groups in the category of {\rm k}-spaces.
\end{prop}

\begin{pf}  Consider the following diagram
$$\xymatrix{
\prod F(X^H,M(G/H))\times \prod F(X^H,M(G/H))\ar[r]^-{\prod +_H}\ar@{->>}[d]_-{p_X\times p_X} &
\prod F(X^H,M(G/H))\ar@{->>}[d]^-{p_X} \\
\EF(X,M)\times \EF(X,M)\ar[r]_-{+}\ar@{->>}[d]_-{\beta_X\times \beta_X} & \EF(X,M)\ar@{->>}[d]^-{\beta_X} \\
\EFG(X,M)\times \EFG(X,M)\ar[r]_-{+} & \,\,\EFG(X,M)\,.}$$
The function on the top is given by the product of the sum on each topological group
$F(X^H,M(G/H))$ and therefore it is continuous.
Since $p_X$ and $\beta_X$ are homomorphisms, both squares commute.  Furthermore,
since $p_X$ and $\beta_X$ are identifications, so are $p_X\times p_X$ and
$\beta_X\times \beta_X$. Therefore $\EF(X,M)$ and $\EFG(X,M)$ are topological
groups.  Finally, since $\OEFG(X,M)$ has the relative k-topology, it is also a topological
group.
\end{pf}

\begin{prop}\label{conthom} Let $f:X\en Y$ be a pointed $G$-map.  Then the homomorphism
 $f_*:\EF(X,M)\en \EF(Y,M)$ is continuous.  Thus the homomorphisms
 $\overline f^G_*:\OEFG(X,M)\en \OEFG(Y,M)$
and $f^G_*:\EFG(X,M)\en \EFG(Y,M)$ are also continuous.
\end{prop}

\begin{pf} The following diagram commutes:
 $$\xymatrix{
\prod_{H\subset G} F(X^H,M(G/H))\ar[r]^-{\prod f^H_*}\ar[d]_-{p_X} &
\prod_{H\subset G} F(Y^H,M(G/H))\ar[d]^{p_Y} \\
\EF(X,M)\ar[r]_-{f_*} & \,\,\EF(Y,M)\,.}$$
Indeed, if $(l_Hx_H)_{H\subset G}\in \prod_{H\subset G} F(X^H,M(G/H))$ is a generator,
then
$$p_Y (\prod f^H_*)((l_Hx_H)_{H\subset G}) = p_Y((l_Hf(x_H))_{H\subset G})= \sum_{H\subset G}
M_*(q_{f(x_H)})(l_H)f(x_H)\,,$$ and
$$f_*p_X((l_Hx_H)_{H\subset G}) = f_*(\sum_{H\subset G}
M_*(q_{x_H})(l_H)x_H) = \sum_{H\subset G}
M_*(\widehat f_{x_H})M_*(q_{x_H})(l_H)f(x_H)\,.$$
Here we denote $q_{H,x_H}$ simply by $q_{x_H}$. Both composites are equal,
since $\widehat f_{x_H}\circ q_{x_H} = q_{f(x_H)}$.

Since $p_X$ is an identification and $\prod f^H_*$ is continuous,
$f_*$ is also continuous.
\end{pf}

The next result shows the homotopy invariance of the functors $\EF$.

\begin{prop}\label{homotinvar} Let $f_0,f_1:X\en Y$ be $G$-homotopic pointed $G$-maps.
Then $f_{0*}, f_{1*}: \EF(X,M)\en \EF(Y,M)$ are $G$-homotopic homomorphisms. Moreover, also
$\overline f^G_{0*}, \overline f^G_{1*}: \OEFG(X,M)\en \OEFG(Y,M)$ and
$f^G_{0*}, f^G_{1*}: \EFG(X,M)\en \EFG(Y,M)$ are homotopic homomorphisms.
\end{prop}

\begin{pf} For each $H\subset G$, take the restriction $f_\nu^H:X^H\en Y^H$, and let
$\ha:X\times I\en Y$ be a $G$-homotopy such that $\ha(x,\nu) = f_\nu(x)$, $\nu = 0,1$,
and denote by $\ha^H:X^H\times I\en Y^H$ the restriction of $\ha$, which is a homotopy between
$f^H_0$ and $f^H_1$.  By \cite{mccord}, there is a homotopy $\widetilde\ha^H:F(X^H,M(G/H))\times I\en
F(Y^H,M(G/H))$ between $(f_0^H)_*$ and $(f_1^H)_*$, which is given by
$\widetilde\ha^H(u,t) = (\ha^H_t)_*(u)$, where $\ha^H_t(x) = \ha^H(x,t)$, $x\in X^H$.

Define a homotopy $R:\prod F(X^H,M(G/H))\times I\en \prod F(Y^H,M(G/H))$ by $R((u_H)_{H\subset G},t) =
(\widetilde\ha^H(u_H,t))_{H\subset G}$.  Define $\ka:\EF(X,M)\times I\en \EF(Y,M)$ by
$\ka(u,t) = \ha_{t*}(u)$ and consider the diagram
$$\xymatrix{
\prod F(X^H,M(G/H))\times I\ar[r]^-{R}\ar@{->>}[d]_{p_X\times \id} &
\prod F(Y^H,M(G/H))\ar@{->>}[d]^-{p_Y} \\
\EF(X,M)\times I\ar[r]_-{\ka} & \,\,\EF(Y,M)\,.}$$
Since both $R$ and $\ka$ are homomorphisms for each fixed $t$, we may
check the commutativity on a generator $((l_Hx_H)_{H\subset G},t)\in \prod F(X^H,M(G/H))\times I$ as follows:
\begin{align*}
p_Y R((l_Hx_H)_{H\subset G},t) & = p_Y((l_H\ha^H(x_H,t))_{H\subset G}) \\
& = \sum_{H\subset G} M_*(q_{\ha^H(x_H,t)})(l_H)\ha^H(x_H,t)\,, \\
\ka(p_X\times \id)((l_Hx_H),t) & = \ka\left(\sum_{H\subset G} M_*(q_{x_H})(l_H)x_H\right) \\
 & = \sum_{H\subset G} M_*(\widehat{\ha_t}_{x_H})M_*(q_{x_H})(l_H)\ha(x_H,t)\,.
\end{align*}
Both expressions are equal, since $\ha^H$ is the restriction of $\ha$, and one clearly has that
$q_{\ha^H(x_H,t)}:G/H\en G/G_{\ha(x_H,t)}$ is equal to the composite $G/H\stackrel{q_{x^H}}{\onto} G/G_{x_H} \stackrel{\widehat{\ha_t}_{x_H}}{\onto}
G/G_{\ha(x_H,t)}$.
Since $p_X\times \id$ is an identification, $\ka$ is continuous an so it is a homotopy as desired.

The homotopy $\ka$ is $G$-equivariant.  Indeed, since $\ha_{t*}$ is $G$-equivariant, one has
$\ka(g\cdot u,t) = \ha_{t*}(g\cdot u) = g\ha_{t*}(u) = g\ka(u,t)$.  Thus we can take the restriction
of $\ka$, $\overline \ka^G:\OEFG(X,M)\times I\en \OEFG(Y,M)$, which is a homotopy between
$\overline f_{0*}^G$ and $\overline f_{1*}^G$, and by the naturality of $\beta$,
we also have a commutative diagram
$$\xymatrix{
\EF(X,M)\times I\ar[r]^-{\ka}\ar[d]_-{\beta_X\times \id} &  \EF(Y,M)\ar[d]^-{\beta_Y} \\
\EFG(X,M)\times I\ar[r]_-{\ka^G} & \,\,\EFG(Y,M)\,,}$$
where $\ka^G(v,t) = \ha^G_{t*}(v)$. Thus $f_{0*}^G, f_{1*}^G:\EFG(X,M)\en \EFG(Y,M)$ are also homotopic.
\end{pf}

We shall need the following $G$-equivariant version of the Whitehead Theorem.

\begin{prop}\label{whitehead}
  Let $Y$ and $Z$ be $G$-spaces and $\fy:Y\en Z$ be a $G$-equivariant weak homotopy equivalence. Let $y_0\in Y$
  and $\fy(y_0)\in Z$ be base points.  Assume that $(Y,y_0)$ and $(Z,\fy(y_0))$ have the $G$-homotopy type of
  pointed $G$-CW-complexes. Then $\fy:(Y,y_0)\en (Z,\fy(y_0))$ is a pointed $G$-homotopy equivalence.
\end{prop}

\begin{pf}
  Consider the square
  $$\xymatrix{
  (Y,y_0)\ar[r]^-{\fy} & (Z,\fy(y_0))\ar[d]^-{\simeq_G} \\
  (C,c_0)\ar[u]^-{\simeq_G}\ar@{-->}[r]_-{\psi} & \,\,(D,d_0)\,,}$$
  where $(C,c_0)$ and $(D,d_0)$ are pointed $G$-CW-complexes and $\psi$ is defined so that the square commutes.
  Since the vertical arrows are $G$-homotopy equivalences, $\psi$ is a $G$-equivariant weak homotopy equivalence.

  Using \cite[II(2.5)]{tD}, one can show that $\psi$ is a pointed $G$-homotopy equivalence.  Therefore $\fy$ is also a
  pointed $G$-homotopy equivalence.
\end{pf}

As a consequence of the previous two results we have the following.

\begin{cor}\label{homotinvar-2}
  Let $Y$ and $Z$ be $G$-spaces and $\fy:Y\en Z$ be a $G$-equivariant weak homotopy equivalence. Let $y_0\in Y$
  and $\fy(y_0)\in Z$ be base points.  Assume that $(Y,y_0)$ and $(Z,\fy(y_0))$ have the $G$-homotopy type of
  pointed $G$-CW-complexes. Then $\fy$ induces a homotopy equivalence of topological groups $\fy^G_*:\EFG(Y,M)\en \EFG(Z,M)$.\qedsymb
\end{cor}

\begin{lem}\label{resolution}
  Let $(X,x_0)$ be a pointed $G$-space of the $G$-homotopy type of a pointed $G$-CW-complex, and let $\sing(X)$ be the
  singular simplicial $G$-set of $X$.  Let $\rho_X:|\sing(X)|\en X$ be given by $\rho([\sigma,s]) = \sigma(s)$, where $\sigma:
  \Delta^q\en X$ and $s\in \Delta^q$.  Then $\rho^G_{X*}:\EFG(|\sing(X)|,M)\en \EFG(X,M)$ is a homotopy equivalence of topological groups.
\end{lem}

\begin{pf}
  By \cite[1.9(e)]{simbredrcm}, we have $|\sing(X)|^H = |\sing(X)^H|$, and clearly $|\sing(X)^H|=|\sing(X^H)|$. Hence
  $\rho_X^H:|\sing(X)|^H\en X^H$ coincides with $\rho_{X^H}:|\sing(X^H)|\en X^H$, which by Milnor's theorem (see \cite{may}) is
  a weak homotopy equivalence. Therefore, $\rho_X$ is a $G$-equivariant weak homotopy equivalence. Hence, by the previous corollary,
  the result follows.
\end{pf}

\begin{defs}\label{EFSIMP}
Assume now that $Q$ is a simplicial pointed $G$-set and consider the simplicial
group $F(Q,M)$ such that for any $n$, $F(Q,M)_n = F(Q_n,M)$. Given a
morphism $\mu:\textbf{m}\en \textbf{n}$ in $\Delta$, we denote by $\mu^Q:Q_n\en Q_m$
the induced pointed $G$-function.  Then we define $\mu^{F(Q,M)} = \mu^Q_*:F(Q,M)_n\en
F(Q,M)_m$.  Then one has two other simplicial groups $\overline F^G(Q,M)$  and $F^G(Q,M)$.
$\overline F^G(Q,M)$  is the simplicial subgroup of $F(Q,M)$ induced by the restricted functor
$\overline F^G(-,M)$ defined above.  $F^G(Q,M)$  is the simplicial quotient group of $F(Q,M)$ induced by
the quotient functor $F^G(-,M)$ defined in page \pageref{functor}.
\end{defs}

\begin{thm}\label{sacabarras} There is a natural isomorphism of topological groups
$$\Psi_Q : |F(Q,M)|\en \EF(|Q|,M)\quad\text{ given by }\quad\Psi_Q([lx,t]) = M_*(q_{x,t})(l)[x,t]\,,$$
where $q_{x,t}:G/G_x\onto G/G_{[x,t]}$ is the canonical projection.
\end{thm}

\begin{pf} In \cite[2.6]{mackey-3} we showed that $\Psi_Q$ is a natural isomorphism
of abelian groups.  Thus we only need to prove that it is a homeomorphism.
First recall \cite[1.9(e)]{simbredrcm} that $|Q^H| = |Q|^H$. Now we are going to see that the following
diagram commutes.
$$\xymatrix{
|\prod_{H\subset G} F(Q^H,M(G/H))|\ar[r]^-{\cong}\ar@{->>}[dd]_-{|p_Q|} &
\prod_{H\subset G} |F(Q^H,M(G/H))|\ar[d]_-{\cong}^-{\prod \psi_{Q^H}} \\
& \prod_{H\subset G} F(|Q^H|,M(G/H))\ar@{->>}[d]^-{p_{|Q|}} \\
|F(Q,M)|\ar[r]_-{\Psi_Q} &  \,\,\EF(|Q|,M)\,.}$$
Take a generator $[(l_H x_H)_{H\subset G},\tau]\in |\prod_{H\subset G} F(Q^H,M(G/H))|$.
Then we can chase it down
and we get $\left[\sum_{H\subset G} M_*(q_{x_H})(l_H)x_H, \tau\right]$, and then to the right
to obtain $\sum_{H\subset G} M_*(q_{x_H,\tau})M_*(\widehat p_{x_H})(l_H)[x_H,\tau]$. If
we first go right and down with $\prod \psi_{Q^H}$, we get $(l_H[x_H,\tau])$, and then
down again with $p_{|Q|}$, we obtain $\sum_{H\subset G} M_*(r_{x_H})(l_H)[x_H,\tau]$, where
$$\xymatrix @-1pc {
G/H\ar[rr]^-{r_{x_H,\tau}}\ar[dr]_-{p_{x_H}} & & G/G_{[x_H,\tau]} \\
& G/G_{x_H}\ar[ur]_-{q_{x_H,\tau}} & }$$
are the canonical projections.  By \cite{may}, the isomorphism on the top
is a homeomorphism, and by \cite{G-homol}, each $\psi_{Q^H}$ is also a homeomorphism. Since
$p_Q$ is a surjective simplicial map, its realization $|p_Q|$ is an identification. And
since $p_{|Q|}$ is also an identification, then $\Psi_Q$ is a homeomorphism.
\end{pf}

\begin{cor}\label{sacabarras-cor} There are natural isomorphisms of topological groups
$$\overline \Psi^G_Q : |\overline F^G(Q,M)|\en \overline \EF^G(|Q|,M)\quad
\text{ and }\quad \Psi^G_Q : |F^G(Q,M)|\en \EF^G(|Q|,M)\,.$$
\end{cor}

\begin{pf} The following is clearly a commutative diagram:
$$\xymatrix{
|\overline F^G(Q,M)|\ar@{^(->}[r]^-{|\iota_Q|}\ar[d]_-{\overline \Psi^G_Q} &
|F(Q,M)|\ar@{->>}[r]^-{|\beta_Q|}\ar[d]^-{\Psi_Q}_-{\cong} & |F^G(Q,M)|\ar[d]^-{\Psi^G_Q} \\
\overline \EF^G(|Q|,M)\ar@{^(->}[r]^-{\iota_{|Q|}} &
\EF(|Q|,M)\ar@{->>}[r]^-{\beta_|Q|} & \EFG(|Q|,M)}$$
and the vertical arrows on the sides are isomorphisms of abelian groups (see \cite[2.6]{mackey-3}).  Since the one in
the middle, by the previous theorem, is a homeomorphism, and the first horizontal arrows are embeddings,
while the second ones are identifications, $\overline \Psi^G_Q$ and $\Psi^G_Q$ are homeomorphisms too.
\end{pf}

The following is the main result of this section.

\begin{thm}\label{mainthm} Let $M$ be a covariant coefficient system for $G$ and $X$ a
  pointed $G$-space of the homotopy type of a $G$-CW-complex. Then the homotopy groups
$$\pi_q(\EFG(X,M))$$
are naturally isomorphic to the (reduced) Bredon-Illman $G$-equivariant homology
groups $\widetilde H_q^G(X;M_*)$ with coefficients in  $M$.
\end{thm}

\begin{pf}  $F^G(\sing(X),M)$ is a simplicial abelian group and hence, by \cite{may},
it is a chain complex with differential $\partial_q^G:F^G(\sing_q(X),M)\en F^G(\sing_{q-1}(X),M)$
given by $\partial^G_q = \sum_{i=0}^q (-1)^i (d_i^{\sing(X)})^G_*$. By \cite[4.5]{mackey-3}
this chain complex is isomorphic to Illman's chain complex $S^G(X,*;M)$ given in \cite[p. 15]{illman},
whose homology is by definition the Bredon-Illman
$G$-equivariant reduced homology of $X$ with coefficients in $M$,
denoted by $\widetilde{H}^G_q(X;M)$.

We shall give an isomorphism
$$H_q(F^G(\sing_*(X),M))\en
\pi_q(\EFG(X,M))\,.$$

To construct the isomorphism, we shall give several isomorphisms as depicted in the
following diagram.
$$\xymatrix @-0.5pc {
H_q(F^G(\sing(X),M))\ar@{<-}[r]^-{i_*}_-{\cong}\ar@{..>}[d] &
\pi_q(F^G(\sing(X),M))\ar[r]^-{\Psi}_-{\cong} &
\pi_q(\sing(|F^G(\sing(X),M)|))\ar[d]^-{\Phi}_-{\cong} \\
\pi_q(\EFG(X,M))\ar@{<-}[r]^-{\cong}_-{(\rho^G_*)_*} &
\pi_q(\EFG(|\sing(X)|,M))\ar@{<-}[r]^-{\cong}_-{\Psi_{\sing(X)*}} &
\pi_q(|F^G(\sing(X),M)|)}$$
By \cite[22.1]{may}, $i_*$ is an isomorphism. In particular, this shows
that every cycle in $\widetilde{H}^G(X;M)$ is represented by a
chain $u$, all of whose faces are zero.
We call this a {\em special chain}.

The homomorphism  $\Psi$, which is given by $\Psi(u)[\tau]=[u,\tau]$, where
$u$ is a special $q$-chain and $\tau\in \Delta^q$, is an
isomorphism, as follows from \cite[16.6]{may}.

In order to define $\Phi$, we must express $\Psi(u)$ as a map
$\gamma:(\Delta[q],\dot{\Delta}[q])\en
(\sing|F^H(\sing(X),M)|,*)$. By the Yoneda lemma, $\gamma$ is the
unique map such that $\gamma(\delta_q) = \Psi(u)$, where
$\delta_q=\id:\overline q\en \overline q$. The homomorphism $\Phi$, defined by
$\Phi[\gamma][f,\tau']=\gamma(f)(\tau')$, for $f\in \Delta[q]_n$ and $\tau'\in \Delta^n$,
is given by the adjunction between the
realization functor and the singular complex functor (see \cite[16.1]{may}).

By Theorem \ref{sacabarras}, $\Psi_{\sing(X)}$ is an isomorphism of topological
groups.

Finally, by Lemma \ref{resolution}, the homomorphism $\rho^G_*$ is a homotopy equivalence.
\end{pf}

\begin{remark} Chasing along the diagram and using the canonical homeomorphism $|\Delta[q]|\en
\Delta^q$ given by $[f,\tau']\mapsto f_\#(\tau')$, one obtains that the isomorphism
$H_q(F^G(\sing(X),M))\en \pi_q(\EFG(X,M))$ is given as follows. It
  maps a homology class $[u]$ represented by a special chain $u=\sum_\alpha
  \gamma^G_{\sigma_\alpha}(u(\sigma_\alpha))$ to the homotopy class $[\overline{u}]$ given by
  $\overline{u}(\tau)=\sum_{\alpha} \gamma^G_{\sigma_\alpha(\tau)}(M_*(p_{\alpha})(u(\sigma_\alpha)))$, where
  $p_\alpha:G/G_{\sigma_\alpha}\en G/G_{\sigma_\alpha(\tau)}$ is the quotient
  function and $\sing(X)/G=\{[\sigma_\alpha]\}$.
\end{remark}

\section[Other topological abelian groups]{\large\sc Other topological abelian groups}\label{sec3}

A different way of topologizing the abelian groups $F(X^\delta,M)$ and $\overline F^G(X^\delta,M) = F^G(X^\delta,M)$,
where $X^\delta$ denotes the underlying set of $X$, is as follows.

\begin{defs}\label{topboletin} Let $X$ be a pointed $G$-space (not necessarily a k-space)
and $M$  be a covariant coefficient system.  Denote by $\sing(X)$ the singular simplicial
pointed $G$-set associated to $X$ and consider the surjective map
$$\pi_X:|F(\sing(X),M)|\onto F(X^\delta,M)\quad\text{ given by }\quad
\pi_X([l\sigma,t]) = M_*(p_{\sigma,t})(l)\sigma(t)\,,$$
where $p_{\sigma,t} : G/G_{\sigma}\onto G/G_{\sigma(t)}$ is the canonical projection.

Give $F(X^\delta,M)$ the identification topology defined by $\pi_X$ and denote the resulting
space by $\F(X,M)$.  Moreover, denote by $\overline \F^G(X,M)$ the group $\overline F^G(X^\delta,M)$
with the relative {\rm k}-topology induced by $\iota_X$ and denote by $\F^G(X,M)$ the group $F^G(X^\delta,M)$ with
the quotient topology induced by $\beta_X$.

Consider the restriction of $\pi_X$
$$\widehat \pi^G_X: |\overline F^G(\sing(X),M)|\onto \overline F^G(X,M)$$
and denote by $\widehat \F^G(X,M)$ the resulting identification space.

We thus have a commutative diagram
$$\xymatrix{
|\overline F^G(\sing(X),M)|\ar@{->>}[d]_-{\widehat \pi^G_X}\ar@{^(->}[rr]^-{|\iota_{\sing(X)}|} & &
|F(\sing(X),M)|\ar@{->>}[d]_-{\pi_X}\ar@{->>}[r]^-{|\beta_{\sing(X)}|} & |F^G(\sing(X),M)|\ar@{->>}[d]^-{\pi^G_X} \\
\widehat \F^G(X,M)\ar[r]_-{\id} & \overline \F^G(X,M)\ar@{^(->}[r]_-{\iota_X} & \F(X,M)\ar@{->>}[r]_-{\beta_X} &
\,\,\F^G(X,M)\,,}$$
where the vertical maps are identifications.  Thus $\id :\widehat \F^G(X,M)\en \overline \F^G(X,M)$ is continuous.
\end{defs}

\begin{remark} The topological groups $\overline F^G(X,M)$ were defined in \cite{bihomol}. \end{remark}

Recall \cite{G-transrami} that a (pointed) $G$-space $X$ is called a {\em  strong
$\rho$-space} if the map $\rho_X:|\sing(X)|\en X$, given by $\rho_X([\sigma,t]) = \sigma(t)$, is a
$G$-retraction.

\begin{prop}\label{huecavsgarigol} If $X$ is a strong $\rho$-space, then $\id:\F(X,M)\en \EF(X,M)$ is a homeomorphism,
and hence the maps $\id:\overline \F^G(X,M)\en \overline \EF^G(X,M)$ and $\id:\F^G(X,M)\en \EFG(X,M)$ are also homeomorphisms.
\end{prop}

\begin{pf} Consider the diagram
$$\xymatrix{
|F(\sing(X),M)|\ar@{->>}[d]_-{\pi_X}\ar[r]^-{\Psi_{\sing(X)}}_-{\cong} &
\EF(|\sing(X)|,M)\ar@{->>}[d]^-{\rho_{X*}} \\
\F(X,M)\ar[r]^-{\cong}_-{\id} & \,\,\EF(X,M)\,.}$$
One easily verifies that the diagram commutes.
The map on the top is a homeomorphism by \ref{sacabarras}, the vertical map on the left is
an identification by definition, and
the vertical map on the right is an identification, since it is a retraction.  Thus the identity
on the bottom is a homeomorphism.
\end{pf}

\begin{cor}\label{garihueca} If $X$ is a strong $\rho$-space,
then the topological groups  $\F^G(X,M)$ and $\EFG(X,M)$ are equal, as well as
the topological groups $\overline \F^G(X,M)$ and $\overline \EF^G(X,M)$.\qedsymb
\end{cor}

\begin{examp} If $K$ is a simplicial $G$-set, then by \cite{G-transrami}, $|K|$ is a strong $\rho$-space.
Hence the topological groups $\F(|K|,M)$ and $\EF(|K|,M)$ are equal, and thus also the
topological groups $\F^G(|K|,M)$ and $\EFG(|K|,M)$ are equal, as well as
the topological groups $\overline \F^G(|K|,M)$ and $\overline \EF^G(|K|,M)$.
\end{examp}

\begin{prop}\label{prop4} If $X$ is a strong $\rho$-space, then $\F^G(X,M) = \EFG(X,M)$.
\end{prop}

\begin{pf} Since $X$ is a strong $\rho$-space, the map $\rho_X:|\sing(X)|\en X$ is a $G$-retraction.
Thus the epimorphism $\rho^G_X:\EFG(|\sing(X)|,M)\onto \EFG(X,M)$
is an identification.  Consider the diagram
$$\xymatrix{
|\FG(\sing(X),M)|\ar[r]^-{\psi^G_M} \ar@{->>}[d]_-{\pi^G_X} & \EFG(|\sing(X)|,M)\ar@{->>}[d]^-{\rho^G_*} \\
\F^G(X,M)\ar[r]_-{\id} & \,\,\EFG(X,M)\,.}$$
By definition of $\pi^G_X$ it commutes and by \ref{sacabarras}, the arrow on the top is an
isomorphism of topological groups.  Since both vertical arrows are identifications, the
result follows.
\end{pf}

There are other topological groups related to the ones previously defined, that were studied in
\cite{bihomol}.

\begin{defs}\label{FG} Let $X$ be a pointed $G$-space (not necessarily a k-space) and $M$ be
a covariant coefficient system for $G$.  Define the topological group (in the category of
k-spaces)
$$F(X,M) = |F(\sing(X),M)|\,,$$
as well as the subgroup $\overline F^G(X,M) = |\overline F^G(\sing(X),M)|$ and the
quotient group $F^G(X,M) = |F^G(\sing(X),M)|$.  One clearly has continuous homomorphisms
$$\overline F^G(X,M)\stackrel{\iota_X}{\hook} F(X,M) \stackrel{\beta_X}{\onto} F^G(X,M)$$
given by $\iota_X = |\iota_{\sing(X)}|$ and $\beta_X = |\beta_{\sing(X)}|$.  The first is
an embedding and the second a quotient map of topological groups.
\end{defs}

\begin{remark} One has identifications of topological groups
$$\overline F^G(X,M)\onto \widehat \F^G(X,M)\,,\quad F(X,M)\onto \F(X,M)\,,\quad
F^G(X,M)\onto \F^G(X,M)\,.$$
We shall see below under which conditions the topological groups $\widehat \F^G(X,M)$ and
$\overline \F^G(X,M)$ are equal.
\end{remark}

\section[Coefficients in $\kmod$]{\large\sc Coefficients in $\kmod$}\label{sec4}

In this section we shall assume that $M:{\mathcal O}(G)\en \kmod$, where $k$ is a field
of characteristic $0$ or a prime $p$ that does not divide the order of $G$,

Recall that if $S$ is a pointed $G$-set, then we have homomorphisms
$$\overline F^G(S,M)\stackrel{\iota_S}{\hook} F(S,M)\stackrel{\beta_S}{\onto} F^G(S,M)\,.$$
We shall study the composite $\alpha_S = \beta_S\circ \iota_S:\overline F^G(S,M)\en F^G(S,M)$.

\begin{prop}\label{prop1} Let $S$ be a pointed $G$-set. Then
 $\alpha_S$ is
a natural isomorphism.
\end{prop}

\begin{pf} Take a generator
$$\gamma^G_x(l) = \sum_{[g]\in G/G_x} M_*(R_{g^{-1}})(l)(gx)\in \overline F^G(S,M)\,.$$
Since $\gamma^G_{gx}M_*(R_{g^{-1}}) = \gamma^G_x$, we have
\begin{align*}
\alpha_S(\gamma^G_x(l)) & = \beta_S \left(\sum_{[g]\in G/G_x} M_*(R_{g^{-1}})(l)(gx)\right) \\
    & = \sum_{[g]\in G/G_x} \beta_S\left(M_*(R_{g^{-1}})(l)(gx)\right) \\
    & = \sum_{[g]\in G/G_x} \gamma^G_{gx}M_*(R_{g^{-1}})(l) \\
    & = \sum_{[g]\in G/G_x} \gamma^G_x(l) \\
    & = [G\!:\!G_x]\gamma^G_x(l) \\
    & = \gamma^G_x([G\!:\!G_x]l)\,.
    \end{align*}
Since $p\not|\,\, |G|$, the indexes $[G\!:\!G_x]$, seen as elements in $k$, are invertible elements, i.e., there exist
the elements $[G\!:\!G_x]^{-1}\in k$.  Since for every generator $\gamma^G_x(l)$, one has
$$\alpha_S(\gamma^G_x([G\!:\!G_x]^{-1}l) = \gamma^G_x(l)\,,$$
$\alpha_S$ is surjective.  To see that it is injective, assume that $u\in \overline F^G(S,M)$
is such that $\alpha_S(u) = 0$.  Hence
$$\alpha_S(u)(x) = [G\!:\!G_x]u(x) = 0 \,\,\forall \, x\,.$$
Thus $[G\!:\!G_x]^{-1}[G\!:\!G_x]u(x) = u(x) = 0$ for all $x$ and so $u = 0$.

The naturality of $\alpha_S$ follows from \eqref{naturality}.
\end{pf}

\begin{prop}\label{prop2} Let $X$ be a pointed $G$-space.  Then \begin{itemize}
\item[{\rm (a)}] $\alpha_X = |\alpha_{\sing(X)}|:\overline F^G(X,M)=|\overline F^G(\sing(X),M)|\en |F^G(\sing(X),M)|=F^G(X,M)$ and
\item[{\rm (b)}] $\alpha_X:\widehat \F^G(X,M)\en \F^G(X,M)$
\end{itemize}
are isomorphisms of topological groups.
\end{prop}

\begin{pf} (a) follows immediately from Proposition \ref{prop1}.  To see (b), consider the
diagram
$$\xymatrix{
|\overline F^G(\sing(X),M)|\ar[r]^-{|\alpha_{\sing(X)}|}_-{\cong}\ar@{->>}[d]_-{\widehat \pi^G_X} &
|F^G(\sing(X),M)^|\ar@{->>}[d]^-{\pi^G_X} \\
\widehat \F^G(X,M)\ar[r]_{\alpha_X}^-{\cong} & \,\,\F^G(X,M)\,.}$$
To see that it is commutative, take a generator $[\gamma^G_\sigma(l),t]\in |\overline F^G(\sing(X),M)|$.
Then
\begin{align*}
\alpha_X\widehat\pi^G_X([\gamma^G_\sigma(l),t]) & =
\alpha_X([G_{\sigma(t)}\!:\!G_\sigma]\gamma^G_{\sigma(t)}M_*(p_\sigma)(l)) \\
& = [G_{\sigma(t)}\!:\!G_\sigma]\gamma^G_{\sigma(t)}M_*(p_\sigma)([G\!:\!G_{\sigma(t)}]l) \\
& = \gamma^G_{\sigma(t)}M_*(p_\sigma)([G\!:\!G_\sigma]l) \\
\\
\pi^G_X|\alpha_{\sing(X)}|([\gamma^G_\sigma(l),t]) & = \pi^G_X([\gamma^G_\sigma([G\!:\!G_\sigma]l),t]) \\
& = \gamma^G_{\sigma(t)}M_*(p_\sigma)([G\!:\!G_\sigma]l)\,.
\end{align*}
Thus we have the assertion.

Since both vertical arrows are identifications, $|\alpha_{\sing(X)}|$ on the top is a homeomorphism
and $\alpha_X$ on the bottom is bijective, then $\alpha_X$ is a homeomorphism too.
\end{pf}

\begin{cor}\label{cor3} The topological groups $\overline \F^G(X,M)$ and $\F^G(X,M)$ are naturally isomorphic.
\end{cor}

\begin{pf}  The isomorphism of topological groups $\alpha_X$ factors as the composite
$$\xymatrix{
\alpha_X:\widehat \F^G(X,M)\ar[r]^-{\iota_X} & \F(X,M)\ar@{->>}[r]^-{\beta_X} & \F^G(X,M)\,.}$$
Therefore $\alpha_X^{-1}\circ \beta_X: \F(X,M)\en \widehat \F^G(X,M)$ is a left inverse for $\iota_X$,
so that $\iota_X$ is an embedding. Hence  $\widehat \F^G(X,M) = \overline \F^G(X,M)\cong \F^G(X,M)$.
\end{pf}

\section[Some applications of the homotopical Bredon-Illman
homology]{\large\sc Some applications of the homotopical Bredon-Illman homology}\label{sec4.5} %\hfill\break

In this section we shall use the topological groups $\EFG(X,M)$ defined in Section \ref{sec2}, together with the $G$-equivariant
weak homotopy equivalence axiom, that we proved in \cite{simbredrcm} using the groups $\FG(X,M)$,
to show that the Bredon-Illman equivariant homology satisfies the (infinite) wedge axiom and to make some calculations.
We start with a general result.

\begin{lem}\label{connpathcomp}
  Let $X$ have the homotopy type of a CW-complex. Then the connected components
  and the path-components of $X$ coincide.
\end{lem}

\begin{pf} Let $C:\Top\en \Set$ be the functor which associates to a space $X$ the set $C(X)$
of its connected components. This is clearly a homotopy functor.

  Now let $\fy:X\en Y$ be a homotopy equivalence, where $Y$ is a CW-complex. Consider the commutative diagram
  $$\xymatrix{
  \pi_0(X)\ar[r]^-{\fy_*}_-{\approx}\ar[d] & \pi_0(Y)\ar[d] \\
  C(X)\ar[r]_-{C(\fy)}^-{\approx} & \,\,C(Y)\,.}$$
  Since $Y$ is locally path-connected, the arrow on the right is a bijection. Hence the
  arrow on the left is also a bijection.
\end{pf}

\begin{defs}\label{weakprod} Let $\Lambda$ be a set of indexes, and let $p(\Lambda)$ be the
set of finite subsets of $\Lambda$. $p(\Lambda)$ is a directed set by inclusion.
Let $\{F_\alpha\mid \alpha\in \Lambda\}$ be a family of pointed k-spaces.
If $A\subset B$ are finite sets of indexes, then we have an inclusion
$\prod_{\alpha\in A}F_\alpha \subset \prod_{\beta\in B}F_\beta$,
defined by putting the base point $*_\beta$ in each factor $F_\beta$, whenever
$\beta\notin A$.
Define
$$\bigoplus_{\alpha\in \Lambda} F_\alpha = \mathop{\colim}_{A\in p(\Lambda)}\prod_{\alpha \in A} F_\alpha\,.$$
Since the category $\ktop$ is closed under colimits, this new space lies in $\ktop$. Notice that when
$\Lambda = \N$ is the set of natural numbers, then the space defined above coincides with the \emph{weak product}
defined in \cite{whitehead}.
\end{defs}

\begin{remark}\label{remweakprod}
  If the spaces $F_\alpha$ in the previous definition are topological abelian groups, then algebraically
  $\bigoplus_{\alpha\in \Lambda} F_\alpha$ is the direct sum of them, and the topology defined therein
  will give this sum a structure of a topological group (see next proposition).  Thus, as topological groups,
  $$\bigoplus_{\alpha\in \Lambda} F_\alpha = \mathop{\colim}_{A\in p(\Lambda)}\bigoplus_{\alpha \in A} F_\alpha\,.$$
\end{remark}

\begin{prop}\label{contprod}
  If for each $\alpha\in \Lambda$, $F_\alpha$ is a T$_1$ topological group, then
  $\bigoplus_{\alpha\in \Lambda} F_\alpha$ is a topological group.
\end{prop}

\begin{pf}
  Since we are taking products in $\ktop$, we consider a compact Hausdorff space $K$
  and any continuous map $f:K\en \bigoplus_{\alpha\in \Lambda} F_\alpha\times
  \bigoplus_{\alpha\in \Lambda} F_\alpha$.  We thus have to prove that the composite
  $$\xymatrix{
  K\ar[r]^-{f} & \bigoplus_{\alpha\in \Lambda} F_\alpha\times
  \bigoplus_{\alpha\in \Lambda} F_\alpha\ar[r]^-{+} & \bigoplus_{\alpha\in \Lambda} F_\alpha}$$
  is continuous.

  To see this, notice that the family $\{\bigoplus_{\alpha\in A} F_\alpha \mid A\in p(\Lambda)\}$
  has the following two properties:
  \begin{itemize}
  \item[(i)] For any $A,B\in p(\Lambda)$ the index $C = A\cap B\in p(\Lambda)$ satisfies
  $$\bigoplus_{\alpha\in A} F_\alpha\cap
  \bigoplus_{\beta\in B} F_\alpha = \bigoplus_{\gamma\in C} F_\alpha\,.$$
  \item[(ii)] For each $A\in p(\Lambda)$, the set $\{B\in p(\Lambda) \mid \bigoplus_{\beta\in B} F_\beta
  \subset \bigoplus_{\alpha\in A} F_\alpha\}$ is finite.
  \end{itemize}
Therefore, by \cite[15.10]{gray}, there are indexes $A_1,\dots,A_m$ and $B_1,\dots, B_n$, such that
$p_1 f(K)\subset \bigoplus_{\alpha\in A_1} F_\alpha\cup\cdots\cup \bigoplus_{\alpha\in A_m} F_\alpha$ and
$p_2 f(K)\subset \bigoplus_{\beta\in B_1} F_\beta\cup\cdots\cup \bigoplus_{\beta\in B_n} F_\beta$,
where $p_1$ and $p_2$ are the projections. Let $C = A_1\cup\cdots \cup A_m\cup B_1\cup\cdots\cup B_n$.
Hence we have the following commutative diagram:
  $$\xymatrix{
  & \bigoplus_{\gamma\in C} F_\gamma \times \bigoplus_{\gamma\in C} F_\gamma \ar[r]^-{+} &
  \bigoplus_{\gamma\in C} F_\gamma\ar@{^(->}[d] \\
  K\ar[r]_-{f}\ar@{-->}[ur] & \bigoplus_{\alpha\in \Lambda} F_\alpha\times
  \bigoplus_{\alpha\in \Lambda} F_\alpha\ar[r]^-{+} & \bigoplus_{\alpha\in \Lambda} F_\alpha}$$
where the sum on the top is continuous because a finite product of topological groups is a topological group, and
the vertical inclusion on the right is continuous.  Hence the composite on the bottom is continuous.
\end{pf}

\begin{prop}\label{finwedge}
  Let $X$ and $Y$ be pointed spaces and $L$ be an abelian group. Then there is an
  isomorphism of topological groups $$F(X\vee Y,L)\cong F(X,L)\times F(Y,L)\,.$$
\end{prop}

\begin{pf}
  Consider the sequences of spaces
  $$\xymatrix{
  X\ar@{^(->}[r]^-{i} & X\vee Y\ar@{->>}[r]^-{q} & Y}\quad
\xymatrix{
  Y\ar@{^(->}[r]^-{j} & X\vee Y\ar@{->>}[r]^-{p} & X}$$
  Then it follows that
  $$\xymatrix{
  F(X,L)\ar[r]^-{i_*} & F(X\vee Y,L)\ar[r]^-{q_*} & F(Y,L)}$$
  is a short exact sequence that splits.  Namely, by functoriality
  one has that $i_*$ is a split monomorphism, $q_*$ is a split epimorphism and
  $q_*\circ i_* = 0$.  Now, if $v\in F(X\vee Y,L)$ is such that $q_*(v)=0$,
  then $q_*(v)(y) = v(x_0,y) = 0$.  Thus $v = i_*(u)$, where $u\in F(X,L)$
  is given by $u(x) = v(x,y_0)$.
  \end{pf}

  Coming back to Definition \ref{weakprod}, we have the following generalization
  to the infinite case of the
  previous proposition.

\begin{prop}\label{infinwedge}
Let $X_\alpha$, $\alpha\in \Lambda$, be a family of pointed spaces.
Then there is an isomorphism of topological groups $\bigoplus_{\alpha\in \Lambda}F(X_\alpha,L)
\cong F(\bigvee_{\alpha\in \Lambda}X_\alpha,L)$.
\end{prop}

\begin{pf}
   For each $A\in p(\Lambda)$, call $\psi_A:\bigoplus_{\alpha\in A}F(X_\alpha,L)\cong
   F(\bigvee_{\alpha\in A}X_\alpha,L)\en F(\bigvee_{\alpha\in \Lambda}X_\alpha,L)$,
   where the isomorphism comes from Proposition \ref{finwedge} and the second map
   is induced by the canonical inclusion.  If $A\subset B$, one easily verifies that
   $\psi_B = \psi_A\circ \psi_{A,B}$, where $\psi_{A,B} : \bigoplus_{\alpha\in A}F(X_\alpha,L)\en
   \bigoplus_{\alpha\in B}F(X_\alpha,L)$ is induced by the inclusion $\bigvee_{\alpha\in A}X_\alpha
   \subset\bigvee_{\alpha\in B}X_\alpha$ modulo the isomorphism of \ref{finwedge}.

   By the universal property of the colimit, the maps (continuous homomorphisms) $\psi_A$ induce a
   continuous homomorphism
   $$\psi:\bigoplus_{\alpha\in \Lambda}F(X_\alpha,L)\en
   F(\bigvee_{\alpha\in A}X_\alpha,L)\,.$$

  In order to see that $\psi$ is an isomorphism of topological groups, we now construct an inverse
  $$\xi : F(\bigvee_{\alpha\in \Lambda}X_\alpha,L)\en
  \bigoplus_{\alpha\in \Lambda}F(X_\alpha,L)$$ as follows. Let $u:\bigvee_{\alpha\in \Lambda}X_\alpha
  \en L$ be an element in $F(\bigvee_{\alpha\in \Lambda}X_\alpha,L)$ and let $u_\alpha :X_\alpha \en L$ be
  the restriction of $u$; that is, $u_\alpha = p_{\alpha *}(u)$, where $p_\alpha :
  \bigvee_{\alpha\in \Lambda} X_\alpha\en X_\alpha$ is the canonical projection. Then $u_\alpha \neq 0$
  only for finitely many values of $\alpha$, i.e., only for $\alpha\in A$, and some $A\in p(\Lambda)$.
  Thus $(u_\alpha)\in \bigoplus_{\alpha\in \Lambda}F(X_\alpha,L)$ and one can define $\xi(u) = (u_\alpha)$.
  The homomorphism $\xi$ is clearly an (algebraic) inverse of $\psi$.  To see that $\xi$ is continuous,
  define the function $\chi:L\times \bigvee_{\alpha\in \Lambda}\en \bigoplus_{\alpha\in \Lambda} F(X_\alpha,L)$
  by $\chi(l,x_\alpha) = lx_\alpha$.  Composing $\chi$ with the identification $\coprod_{\alpha\in \Lambda}
  L\times X_\alpha\onto L\times \bigvee_{\alpha\in \Lambda} X_\alpha$ and then restricting to each
  $L\times X_\alpha$, we obtain the composite
  $$L\times X_\alpha\stackrel{\mu_1}{\en} F(X_\alpha,L)\hook \bigoplus_{\alpha\in \Lambda} F(X_\alpha,L)$$
  which is continuous.  Therefore $\chi$ is continuous.  Now consider the commutative diagram
  $$\xymatrix{
  (L\times \bigvee_{\alpha\in \Lambda} X_\alpha)^k\ar[r]^-{\chi^k}\ar@{->>}[d] &
  (\bigoplus_{\alpha\in \Lambda} F(X_\alpha,L))^k\ar[d]^-{+} \\
  F(\bigvee_{\alpha\in \Lambda}X_\alpha,L)\ar[r]_-{\xi} & \,\,\bigoplus_{\alpha\in \Lambda} F(X_\alpha,L)\,.}$$
  Since the vertical map on the left is an identification and the maps on the top and on the
  right are continuous ($+$ by \ref{contprod}), $\xi$ is continuous.
  \end{pf}

Now we can use the groups $\EFG(X,M)$ to prove that the Bredon-Illman
$G$-equiv\-ariant  homology theory $\widetilde H^G_*(-;M)$ satisfies the
wedge axiom.  We need a lemma whose proof is not difficult.

\begin{lem}\label{colimidentif}
  Let $\{Y_\alpha; i_{\beta,\alpha}:Y_\alpha\hook Y_\beta\}$ and
  $\{\overline Y_\alpha; \overline i_{\beta,\alpha}:\overline Y_\alpha\hook \overline Y_\beta\}$
  be diagrams in $\ktop$ such that each $i_{\beta,\alpha}$ and $\overline i_{\beta,\alpha}$ are closed embeddings,
  and let  $\{q_\alpha : Y_\alpha \onto \overline Y_\alpha\}$ be a family of
  identifications such that $q_\beta\circ i_{\beta,\alpha} = \overline i_{\beta,\alpha}\circ q_\alpha$.
  Then the map $q:\colim Y_\alpha\onto \colim \overline Y_\alpha$ induced in the colimits in $\ktop$ is
  an identification.\qedsymb
\end{lem}

The next lemma is a consequence of the previous one.

\begin{lem}  Let $X_\alpha$ be a pointed $G$-space for each $\alpha\in \Lambda$.
  Then the map
  $$\xymatrix{\bigoplus_{\alpha\in \Lambda} \prod_{H\subset G} F(X_\alpha^H;M(G/H))\ar[r]^-{\oplus p_{X_\alpha}} &
  \bigoplus_{\alpha\in \Lambda} \EFG(X_\alpha,M)}$$
  is an identification.\qedsymb
\end{lem}

Now, as an application of our groups $\EFG(X,M)$, we have the next results. First we have that
the Bredon-Illman homology satisfies the wedge-axiom.

\begin{prop}\label{wedge}
  Let $X_\alpha$, $\alpha\in \Lambda$, be an arbitrary family of pointed $G$-spaces.  Then there is an isomorphism
  $$\EFG(\bigvee_{\alpha\in \Lambda} X_\alpha,M)\cong \bigoplus_{\alpha\in \Lambda} \EFG(X_\alpha,M)\,.$$
  \end{prop}

  \begin{pf}
    By Proposition \ref{infinwedge}, the inclusions $i_\alpha : X_\alpha\hook \bigvee_\alpha X_\alpha$ induce
    an isomorphism $\fy$ of topological groups by
    $$\xymatrix{
    \bigoplus_{\alpha\in \Lambda} \prod_{H\subset G} F(X_\alpha^H;M(G/H))\ar@{-->}[dr]^-{\fy}\ar[d]_-{\cong} & \\
    \prod_{H\subset G}\bigoplus_{\alpha\in \Lambda}F(X_\alpha^H,M(G/H))\ar[r]_-{\prod \psi_H}^-{\cong} &
    \prod_{H\subset G}F(\bigvee_{\alpha\in \Lambda} X_\alpha^H,M(G/H))}$$
    Similarly, we can define an isomorphism of abelian groups
    $$\xymatrix{\psi^G:\bigoplus_{\alpha\in \Lambda} \EFG(X_\alpha;M)\ar[r] &
    \EFG(\bigvee_{\alpha\in \Lambda} X_\alpha,M)\,.}$$
    They fit into a commutative
    diagram
    $$\xymatrix{
    \bigoplus_\alpha \prod_{H\subset G} F\left(X_\alpha^H;M(G/H)\right)\ar[r]^-{\fy}\ar@{->>}[d]_-{\oplus p_{X_\alpha}} &
    \prod_{H\subset G} F\left(\bigvee_\alpha X_\alpha^H,M(G/H)\right)\ar@{->>}[d]^-{p_{\vee X_\alpha}} \\
    \bigoplus_\alpha \EFG\left(X_\alpha,M\right)\ar[r]_-{\psi^G} & \,\,\EFG\left(\bigvee_\alpha X_\alpha,M\right)\,,}$$
    where the vertical arrows are identifications (the left one by the previous lemma).
    Since the top arrow is a homeomorphism, then
    the bottom arrow is a homeomorphism too.
      \end{pf}

\begin{prop}\label{prewedgeaxiom}
  Let $X_\alpha$, $\alpha\in \Lambda$, be an arbitrary
  family of pointed $G$-spaces having the homotopy type of $G$-CW-complexes.  Then there is an isomorphism
  $$\widetilde H^G_*(\bigvee_\alpha X_\alpha;M)\cong \bigoplus_\alpha \widetilde H^G_*(X_\alpha;M)\,.$$
\end{prop}

\begin{pf} Under conditions (i) and (ii) given in the proof of \ref{contprod},
that are satisfied by the family $\{\bigoplus_{\alpha\in A} \EFG(\bigvee_\alpha X_\alpha,M) \mid A\in p(\Lambda)\}$,
it is proved in \cite{gray} that the homotopy groups commute with the colimit. Hence
\begin{align*}
  \widetilde H^G_*(\bigvee_{\alpha\in \Lambda} X_\alpha;M) &
  \cong \pi_q(\EFG(\bigvee_{\alpha\in \Lambda} X_\alpha;M)) & \text{(by \ref{mainthm})}\\
  &\cong   \pi_q(\bigoplus_{\alpha\in \Lambda} \EFG(X_\alpha,M)) & \text{(by \ref{wedge})}\\
  &\cong \pi_q(\colim_A (\bigoplus_{\alpha\in A} \EFG(X_\alpha,M))) & \text{(by \ref{remweakprod})}\\
  & \cong \colim_A (\pi_q(\bigoplus_{\alpha\in A} \EFG(X_\alpha,M))) & \text{(by \cite[15.9]{gray})}\\
& \cong \colim_A (\bigoplus_{\alpha\in A}
\pi_q(\EFG(X_\alpha,M))) & \text{(since $A$ is finite)}\\
  & \cong \bigoplus_{\alpha\in \Lambda}\pi_q(\EFG(X_\alpha,M))) & \text{(by \ref{remweakprod})}\\
  & \cong \bigoplus_{\alpha\in \Lambda} \widetilde H^G_*(X_\alpha;M)& \text{(by \ref{mainthm}).}
\end{align*}
\end{pf}

In order to prove the wedge axiom in full generality, we need the following lemmas.

\begin{lem}\label{wedge-fin}
  Let $X$ and $Y$, and $X'$ and $Y'$ be well-pointed spaces and let $\fy:X'\en X$ and $\psi:Y'\en Y$
  be weak homotopy equivalences (mapping base points to base points).  Then $\fy\vee\psi: X'\vee Y'\en
  X\vee Y$ is a weak homotopy equivalence.
\end{lem}

\begin{pf}
  Consider the double attaching spaces $X\cup I\cup Y$ and $X'\cup I\cup Y'$, where the base point
  $x_0\in X$ is identified with $0\in I$, the base point $y_0\in Y$ is identified with $1\in I$,
  and similarly with the other union. Since the spaces are well
  pointed, the quotient maps $q:X\cup I\cup Y\onto X\vee Y$ and $q':X'\cup I\cup Y'\onto X'\vee Y'$ that
  collapse $I$ to the common base point are homotopy equivalences. The pairs $(X\cup [0,1),(0,1]\cup Y)$
  and $(X'\cup [0,1),(0,1]\cup Y')$ are excisive in $X\cup I\cup Y$ and $X'\cup I\cup Y$, respectively.
  Hence, by \cite[16.24]{gray}, the map $\fy\cup \id_I\cup \psi:X\cup I\cup Y\en X'\cup I\cup Y'$ is a weak
  homotopy equivalence.  The following is clearly a commutative diagram
  $$\xymatrix{
            X'\cup I\cup Y'\ar[rr]^-{\fy\cup \id_I\cup \psi}\ar[d]_-{q'}^-{\simeq} & &
            X\cup I\cup Y\ar[d]^-{q}_-{\simeq} \\
            X'\vee Y'\ar[rr]_-{\fy\vee\psi} & &\,\,X\vee Y\,.}$$
  Therefore, $\fy\vee\psi$ is a weak homotopy equivalence.
\end{pf}

As a consequence, we obtain the following result.

\begin{lem}\label{wedge-infin}
  Let $X_\alpha$ and $X'_\alpha$, $\alpha\in \Lambda$ be an arbitrary family of well-pointed $G$-spaces,
  and let $\fy_\alpha:X'_\alpha\en X_\alpha$ be an equivariant weak homotopy equivalences (mapping base points to base points).
  Then $\vee_{\alpha\in \Lambda}\fy_\alpha: \bigvee_{\alpha\in \Lambda}X'_\alpha\en \bigvee_{\alpha\in \Lambda}X_\alpha$
  is an equivariant weak homotopy equivalence.
\end{lem}

\begin{pf}
  First notice that since the base points are fixed under the $G$-action, $(\bigvee_{\alpha\in \Lambda}X_\alpha)^H =
  \bigvee_{\alpha\in \Lambda}X_\alpha^H$, and similarly  $(\bigvee_{\alpha\in \Lambda}X'_\alpha)^H =
  \bigvee_{\alpha\in \Lambda}X_\alpha^{\prime H}$.

  Let $p(\Lambda)$ be the set of finite subsets of $\Lambda$. Then, by \ref{wedge-fin}, for every $A\in p(\Lambda)$,
  $\vee_{\alpha\in A}\fy_\alpha^H: \bigvee_{\alpha\in A}X^{\prime H}_\alpha\en \bigvee_{\alpha\in A}X^H_\alpha$ is a
  weak homotopy equivalence and induces isomorphisms between all the homotopy groups. Thus again, as in the proof
  of \ref{prewedgeaxiom} and using \cite[15.9]{gray}, $\vee_{\alpha\in \Lambda}\fy^H_\alpha:
  \bigvee_{\alpha\in \Lambda}X^{\prime H}_\alpha\en
  \bigvee_{\alpha\in \Lambda}X^H_\alpha$ induces isomorphisms in all homotopy groups and hence it is a weak homotopy
  equivalence, thus $\vee_{\alpha\in \Lambda}\fy_\alpha: \bigvee_{\alpha\in \Lambda}X'_\alpha\en \bigvee_{\alpha\in \Lambda}X_\alpha$
  is an equivariant weak homotopy equivalence.
\end{pf}

Now we can prove the general wedge axiom for the Bredon-Illman $G$-equivariant homology.

\begin{thm}\label{infwedgehomol}
  Let $X_\alpha$, $\alpha\in \Lambda$, be an arbitrary family of well-pointed $G$-spaces.
  Then there is an isomorphism $\bigoplus_{\alpha\in \Lambda}\widetilde H^G_q(X_\alpha;M)\cong
  \widetilde H^G_q(\bigvee_{\alpha\in \Lambda} X_\alpha;M)$ induced by the canonical inclusions
  $X_\alpha\hook \bigvee X_\alpha$.
\end{thm}

\begin{pf} For each $\alpha\in \Lambda$ there is an equivariant weak homotopy equivalence $\fy_\alpha:\widetilde X_\alpha
\en X_\alpha$, where $\widetilde X_\alpha$ is a (pointed) $G$-CW-complex (for instance $\widetilde X_\alpha = |\sing(X_\alpha)|$).

By \ref{wedge-infin}, $\fy=\vee_\alpha \fy_\alpha : \bigvee_\alpha \widetilde X_\alpha
\en \bigvee_\alpha X_\alpha$ is an equivariant weak homotopy equivalence. Hence, by \cite[1.19]{simbredrcm},
$\fy$ induces an isomorphism $\fy_*:H^G_q(\bigvee_\alpha \widetilde X_\alpha;M)\en H^G_q(\bigvee_\alpha X_\alpha)$.
Recall that
$$H^G_q(\bigvee_\alpha \widetilde X_\alpha;M) = \widetilde H^G_q((\bigvee_\alpha \widetilde X_\alpha)^+;M)\quad\text{and}
\quad H^G_q(\bigvee_\alpha X_\alpha)= \widetilde H^G_q((\bigvee_\alpha X_\alpha)^+;M)\,.$$

On the other hand, the cofiber sequences
$$\s^0\hook (\bigvee_\alpha \widetilde X_\alpha)^+\onto \bigvee_\alpha \widetilde X_\alpha\quad\text{and}
\quad\s^0\hook (\bigvee_\alpha X_\alpha)^+\onto \bigvee_\alpha X_\alpha$$
induce short exact sequences that fit into the diagram
$$\xymatrix{
0\ar[r] & \widetilde H^G_0(\s^0;M)\ar[r]\ar@{=}[d] &
\widetilde H^G_0((\bigvee_\alpha \widetilde X_\alpha)^+;M)\ar[r]\ar[d]^-{\fy^+_*}_-{\cong} &
\widetilde H^G_0(\bigvee_\alpha \widetilde X_\alpha;M)\ar[r]\ar[d]^-{\fy_*} & 0 \\
0\ar[r] & \widetilde H^G_0(\s^0;M)\ar[r] &
\widetilde H^G_0((\bigvee_\alpha X_\alpha)^+;M)\ar[r] &
\widetilde H^G_0(\bigvee_\alpha X_\alpha;M)\ar[r] & 0}$$
Thus, by the five-lemma, we have an isomorphism
$$\fy_*:\widetilde H^G_0(\bigvee_\alpha \widetilde X_\alpha;M)\en
\widetilde H^G_0(\bigvee_\alpha X_\alpha;M)\,.$$
On the other hand, since $\widetilde H^G_q(\s^0;M) = 0$ if $q>0$,
there are isomorphisms
$$\widetilde H^G_q((\bigvee_\alpha \widetilde X_\alpha)^+;M)\cong
\widetilde H^G_q(\bigvee_\alpha X_\alpha;M)$$
$$\widetilde H^G_q((\bigvee_\alpha \widetilde X_\alpha)^+;M)\cong \widetilde H^G_q(\bigvee_\alpha X_\alpha;M)$$
and thus also an isomorphism
$$\fy_*:\widetilde H^G_q(\bigvee_\alpha \widetilde X_\alpha;M)\en
\widetilde H^G_q(\bigvee_\alpha X_\alpha;M)\,.$$
Similarly, for every $\alpha\in \Lambda$, there are isomorphisms
$$\fy_{\alpha*}:\widetilde H^G_q(\widetilde X_\alpha;M)\en \widetilde H^G_q(X_\alpha;M)\,.$$

Since the $G$-spaces $\widetilde X_\alpha$ are $G$-CW-complexes, by \ref{prewedgeaxiom} there is an isomorphism
$$\bigoplus_{\alpha\in \Lambda}\widetilde H^G_q(\widetilde X_\alpha;M)\cong \widetilde H^G_0(\bigvee_\alpha \widetilde X_\alpha;M)$$
for all $q\ge 0$, and by the isomorphisms above it follows that there is an isomorphism
$$\bigoplus_{\alpha\in \Lambda}\widetilde H^G_q(X_\alpha;M)\cong \widetilde H^G_q(\bigvee_\alpha X_\alpha;M)$$
as desired.
\end{pf}

Next we make some calculations using the homotopical approach to Bredon-Illman homology.

\begin{prop}\label{hfixconnred}
Let $X$ be a pointed $G$-space of the homotopy type of a $G$-CW-complex. Assume that
$X^H$ is connected for each $H\subset G$.  Then $$\widetilde H^G_0(X;M) = 0\,.$$
\end{prop}

\begin{pf}
By Lemma \ref{connpathcomp}, $C(X^H) = \pi_0(X^H)$ for each $H\subset G$. Hence
$$\pi_0(F(X^H,M(G/H))) \cong \widetilde H_0(X^H;M(G/H) = 0\,.$$
Thus the topological groups $F(X^H,M(G/H))$ are path-connected and so is their
product.  Since $\EFG(X,M)$ is a quotient of this product, its also path-connected.
Hence $0 = \pi_0(\EFG(X,M)) \cong \widetilde H^G_0(X;M)$.
\end{pf}

\begin{prop}\label{hfixconn} Let $X$ be a $G$-space of the homotopy type of a $G$-CW-complex. Assume that
$X^H$ is connected for each $H\subset G$.  Consider the family $\ha$ of all $H\subset G$ such that
$X^H\ne \emptyset$. Then $H^G_0(X;M)$ is a quotient of $\bigoplus_{H\in \ha} M(G/H)$.  Furthermore,
if $X$ has a fixed point, then $H^G_0(X;M)\cong M(G/G)$.
\end{prop}

\begin{pf}
  For each $H\in \ha$, there is an isomorphism
  $$\pi_0(F((X^H)^+,M(G/H))) \cong H_0(X^H;M(G/H))\cong M(G/H)\,.$$
  Hence one has an epimorphism
  $$\bigoplus_{H\in \ha} M(G/H)\cong \pi_0 (\prod_{H\in \ha} F((X^H)^+,M(G/H)))\onto \pi_0(\EFG(X^+,M)) \cong H^G_0(X;M)\,.$$

  Furthermore, if $X$ has a fixed point $x_0$, then $X^+ = X\vee \s^0$, where we take $x_0$ as base point of $X$.
  By the exactness axiom, one always has $\widetilde H^G_*(X\vee Y;M)\cong
  \widetilde H^G_*(X;M)\oplus \widetilde H^G_*(Y;M)$. Then, in particular,
  $\widetilde H^G_0(X^+;M) \cong \widetilde H^G_0(X;M)\oplus
  \widetilde H^G_0(\s^0;M)$.  Since by the previous proposition
  $\widetilde H^G_0(X;M) = 0$, one obtains
  $H^G_0(X;M) \cong \widetilde H^G_0(\s^0;M) \cong M(G/G)$.
\end{pf}

\section[Dold-Thom topological groups and the transfer]{\large\sc Dold-Thom
topological groups and the transfer}\label{sec5}

In this section we summarize the properties of the topological groups defined
and explain which of them admit a transfer, either for ordinary or for ramified
covering $G$-maps.  We start with the following concept.

\begin{defs}\label{DTTGF}  Let $M$ be a covariant coefficient system for $G$.
We shall call a functor $\ga(-,M):\Gtop\en \Topab$
a \emph{Dold-Thom topological group functor with coefficients in $M$ for a
subcategory ${\mathcal T}\subset \Gtop$} if there is a natural isomorphism
$$\fy_X:\pi_*(\ga(X,M))\en \widetilde H^G_*(X;M)\,,$$
where the right-hand side denotes the reduced $G$-equivariant Bredon-Illman homology
groups of $X$ with coefficients in $M$, and $X$ is an object of ${\mathcal T}$.
\end{defs}

\begin{examps} The following are examples of Dold-Thom topological group functors
with coefficients in a coefficient system $M$:
\begin{itemize}
\item[1.] The groups $F^G(X,M) = |F^G(\sing(X),M)|$, for any $G$-space $X$, defined
in \cite{mackey-3}.
\item[2.] The groups $\F^G(X,M)$ defined in Section 5 of \cite{bihomol}, when $M$ is a homological Mackey functor and $X$ has the homotopy type of a $G$-CW-complex.
\item[3.] The groups $\EFG(X,M)$ as shown in \ref{mainthm}.
\item[4.] The groups $|F^G(Q,M)|$ for $X = |Q|$ and $Q$ a simplicial $G$-set, in particular, the simplicial $G$-set associated to a simplicial $G$-complex as shown in \cite{simbredrcm}.
\item[5.] The groups ${\rm AG}(X)$ and ${\rm AG}(X;m)$ defined by Dold and Thom \cite{DT} in the nonequivariant case,
where the coefficients are a cyclic group ($\Z$ and $\Z/m\Z$, respectively) and $X$ is a (countable) simplicial complex. Here G does not stand for any group.
\item[6.] The groups $B(L,X)$ defined by McCord \cite{mccord} in the nonequivariant case and coefficients in an abelian group $L$, where $X$ is a weak Hausdorff k-space of the homotopy type of a CW-complex.
\item[7.] The groups ${\rm AG}(X)$ defined by Lima-Filho \cite{Lima} for a $G$-CW-complex $X$ and coefficients in $\Z$.
\item[8.] The groups $L\otimes X$ defined by dos Santos \cite{santos} for a pointed $G$-CW-complex $X$ and coefficients in a $G$-module $L$.
\item[9.] The groups $\ga X\otimes_{G\EF} M$ defined in \cite{nie} for a pointed $G$-CW-complex $X$.
\end{itemize}
And if $M$ takes values in $k$-$\Mod$, where $k$ is a field of characteristic $0$ or a prime $p$ that does not divide $|G|$:
\begin{itemize}
\item[10.] The groups $\overline F^G(X,M)$, by \ref{prop2}(a) and 1.
\item[11.] The groups $\overline \F^G(X,M)$, if $X$ has the homotopy type of a $G$-CW-complex, by \ref{cor3} and 2.
\item[12.] The groups $\OEFG(X,M)$ for a strong $\rho$-space $X$, by \ref{prop2}(b), \ref{cor3}, \ref{prop4}, and 3.
\end{itemize}
\end{examps}

\begin{defs}\label{transfer}  Let $M$ be a Mackey functor for the finite group $G$ and
$\ga(-,M)$ be a Dold-Thom topological group functor with coefficients in (the covariant part of) $M$.
Let $p:E\en X$ be either an $n$-fold $G$-equivariant ordinary covering map (see \cite{mackey-3})
or an $n$-fold $G$-equivariant ramified covering map (see \cite{G-transrami}).  By a \emph{transfer}
for $p$ in $\ga(-,M)$ we understand a continuous homomorphism
$$t^G_p:\ga(X^+,M)\en \ga(E^+,M)\,,$$
which satisfies the following conditions:
\begin{itemize}
\item[(a)] Pullback: If $f:Y\en X$ is continuous and we take the pullback diagram
$$\xymatrix{
E'\ar[d]_-{q}\ar[r]^-{\widetilde f} & E\ar[d]^-{p} \\
Y\ar[r]_-{f} & \,\,X\,,}$$
then $t^G_p\circ f^G_* = \widetilde f^G_*\circ t^G_q : \ga(Y^+,M)\en \ga(E^+,M)$.
\item[(b)] Normalization: If $p=\id_X:X\en X$, then $t^G_{\id_X} = \id:\ga(X^+,M)\en \ga(X^+,M)$.
\item[(c)] Functoriality: If $p:E\en X$ and $q:X\en Y$ are $G$-equivariant
ordinary, resp. ramified covering maps, then
$$t^G_{q\circ p} = t^G_p\circ t^G_q:\ga(Y^+,M)\en \ga(E^+,M)\,.$$
\item[(d)] If $M$ is homological, then the composite $p^G_*\circ t^G_p:\ga(X^+,M)\en \ga(X^+,M)$
is multiplication by $n$.
\end{itemize}
\end{defs}

\begin{examps} The following Dold-Thom topological group functors have transfers for $p$:
\begin{itemize}
\item[1.] If $p:E\en X$ is an $n$-fold $G$-equivariant ordinary covering map and $M$
is any Mackey functor, then there is a transfer $t^G_p:F^G(X^+,M)\en F^G(E^+,M)$, as shown in
\cite{mackey-3}.
\item[2.] If $p:E\en X$ is an $n$-fold $G$-equivariant ramified covering map between strong $\rho$-spaces of the homotopy type of $G$-CW-complexes, and $M$ is homological, then there is a transfer $t^G_p:\F^G(X^+,M)\en \F^G(E^+,M)$, as shown in
\cite{G-transrami}.
\item[3.] If $p:K\en Q$ is a special $G$-equivariant simplicial ramified covering map and $M$ is any Mackey functor,
then there is a transfer $|t^G_p|:|F^G(Q^+,M)|\en |F^G(K^+,M)|$, as shown in
\cite{simbredrcm}.
\item[4.] If $p:E\en X$ is an $n$-fold $G$-equivariant ramified covering map between strong $\rho$-spaces of the homotopy type of $G$-CW-complexes, and $M$ is homological, then there is a transfer $t^G_p:\EFG(X^+,M)\en \EFG(E^+,M)$, by
\ref{huecavsgarigol} and Example 2.
\end{itemize}
\end{examps}

\begin{remark}
  Let $p:E\en X$ be an $n$-fold $G$-equivariant covering map.  The restrictions $p^H:E^H\en X^H$ are not, in general, covering maps.
  Thus there are no transfers $F(X^{H+},M(G/H))\en F(E^{H+},M(G/H))$.
  Since the topology of the groups $\EFG(X^+,M)$ and $\EFG(E^+,M)$ is given in terms of that of the groups $F(X^{H+},M(G/H))$ and $F(E^{H+},M(G/H))$, it does not seem possible to prove the continuity  of $t^G_p$.   And even if the transfers $t_{p^H}$ exist, they do not commute with the identifications. However, if as stated in Example 4 above, the spaces are $\rho$-spaces, then the groups $\EFG(X^+,M)$ and $\EFG(E^+,M)$ coincide with the groups $\F^G(X^+,M)$ and $\F^G(E^+,M)$, as shown in \ref{garihueca}, and one can show that the transfer is continuous.
\end{remark}

Now we study the homotopy type of the Dold-Thom topological groups. First we have the following general result.

\begin{thm}\label{mooregeneral} Let $A$ be a locally connected
topological abelian group in the category
$\ktop$ that has the homotopy type of a CW-complex. Then there is a homotopy equivalence
$$A\simeq \bigoplus_{q\ge 0} K(\pi_q(A),q)\,,$$
where $K(\pi_q(A),q)$ denotes the corresponding Eilenberg-Mac Lane space.
\end{thm}

\begin{pf} Since translation by $a_0$, $a_0\in A$, is a homeomorphism in the k-topology,
the connected component $A_0$ of $0\in A$ is a closed subgroup of $A$.

By Lemma \ref{connpathcomp}, the connected components of $A$ coincide with the path-com\-ponents.
Thus $A$ is the topological sum of all its path-components, and they are closed and open, because
$A$ is locally connected.
Since all path-components of $A$ are homeomorphic (via translation) to
$A_0$, we have a homeomorphism
\setcounter{thm}{7}
\begin{equation}\label{eq6}
A\cong \pi_0(A)\times A_0\,.\end{equation}
Since $A$ has the homotopy type of a CW-complex, so does $A_0$. Consider
$$A_0\stackrel{i}{\hook} F(A_0,\Z)\stackrel{\nu}{\onto} A_0\,,$$
given by $i(a) = 1a$ and $\nu(u) = \sum_{a\in A_0} u(a)a$.  Then $i$ and $\nu$
are clearly continuous, and the composite
$\nu\circ i$ is equal to the identity $\id_{A_0}$.   Applying the functor
$\pi_q$ we have the following
$$\xymatrix{
\pi_q(A_0)\ar@/^1.2pc/[rr]^-{1}\ar[r]_-{i_\#}\ar[dr]_-{h_q} & \pi_q(F(A_0,\Z))\ar[d]^-{\cong}\ar[r]_-{\nu_\#} &
\pi_q(A_0)  \\
& \widetilde H_q(A_0;\Z)\ar[ru]_-{\lambda_q} & }$$
where $\lambda_q$ is defined so that the triangle commutes, and thus we obtain a
left inverse for the Hurewicz homomorphism.

Now consider a homotopy equivalence $A_0\stackrel{\fy}{\en} C_0$, where $C_0$ is a CW-complex,
and the diagram
$$\xymatrix{
\pi_q(A_0)\ar[r]^-{h_q}\ar[d]_-{\fy_\#}^-{\cong} & \widetilde H_q(A_0)\ar[r]^-{\lambda_q}
\ar[d]_-{\fy_\#}^-{\cong} & \pi_q(A_0)\ar[d]_-{\fy_\#}^-{\cong} \\
\pi_q(C_0)\ar[r]_-{h_q} & \widetilde H_q(C_0)\ar@{-->}[r]_-{\alpha_q} & \,\,\pi_q(C_0)\,,}$$
where $\alpha_q$ is so that the diagram commutes.  Hence $\alpha_q$ is a left inverse to Hurewicz
too.  Since $C_0$ is a connected CW-complex, by Moore's theorem (see \cite[IX(1.9)]{whitehead},
we have $C_0 \simeq \bigoplus_{q\ge 1}K(\pi_q(C_0),q)$ (see \ref{weakprod}).  Hence, using \eqref{eq6},
we get
$$A\cong \pi_0(A)\times A_0\simeq \pi_0(A)\times \bigoplus_{q\ge 1}K(\pi_q(C_0),q)\approx
\bigoplus_{q\ge 0}K(\pi_q(A_0),q)\,.$$
\end{pf}

\begin{remark}
  This result was proved by John Moore in \cite{moore-2} for the case of connected simplicial abelian groups and furthermore by Dold and Thom \cite{DT-2} for the case of connected finite polyhedral abelian groups.
\end{remark}

As a consequence of Theorem \ref{mooregeneral} we have the following.

\begin{thm}\label{uniqueness}
  Let $\ga(-,M)$ is a Dold-Thom topological group functor.  If $X$ is a pointed $G$-space
  such that $\ga(X,M)$ is locally connected and has the homotopy type of a CW-complex,
  then
  $$\ga(X,M)\simeq \bigoplus_{q\ge 0} K(\widetilde H^G_q(X;M),q)\,.$$
  Hence these groups are unique up to homotopy. \qedsymb
\end{thm}

\begin{prop}\begin{itemize}
\item[{\rm (a)}] The topological groups $F^G(X,M)$ are CW-complexes for any pointed $G$-space $X$.
\item[{\rm (b)}] If $X$ is a pointed $G$-CW-complex, then $\EFG(X,M)$ is locally connected and
has the homotopy type of a CW-complex.
\item[{\rm (c)}] If $X$ is a pointed $G$-simplicial complex
or a finite-dimensional countable locally finite $G$-CW-complex, then $\F^G(X,M)$ is locally connected and
has the homotopy type of a CW-complex.
\end{itemize}
We conclude that in all these cases the topological groups have the homotopy
type of $\bigoplus_{q\ge 0} K(\widetilde H^G_q(X;M),q)$.
\end{prop}

\begin{pf} (a) follows from the fact that $F^G(X,M)$ is the geometric realization
$|F^G(\sing(X),M)|$.

(b) First notice that the property of a space being locally connected is inherited by quotient spaces.
Hence, if $X$ is locally connected, so is also $\prod_{n\ge 1}(L\times X)^n$ for any abelian group $L$,
and given the quotient map $\prod_{n\ge 1}(L\times X)^n\onto F(X,L)$, the topological group $F(X,L)$
is locally connected.

Now, if $X$ is a $G$-CW-complex, then $X^H$ is locally connected for every $H\subset G$.
Hence each topological group $F(X^H,M(G/H))$ is locally connected and since by definition
there is a quotient map $\prod_{H\subset G}F(X^H,M(G/H))\onto \EFG(X,M)$,
the topological group $\EFG(X,M)$ is also locally connected. Furthermore there is a
$G$-homotopy equivalence $\rho_X:|\sing(X)|\en X$, which by the homotopy invariance \ref{homotinvar}
induces a homotopy equivalence $\rho_{X*}^G:\EFG(|\sing(X)|,M)\en \EFG(X,M)$.   But
by \ref{sacabarras-cor}, there is an isomorphism of topological groups $\EFG(|\sing(X)|,M)\cong
|\FG(\sing(X),M)|$, thus the first is a CW-complex, and hence $\EFG(X,M)$ has the homotopy
type of a CW-complex.

(c) Actually we only need $X$ to be a pointed $G$-CW-complex which is also a strong $\rho$-space, then by \ref{garihueca},
$\F^G(X,M) = \EFG(X,M)$ and the result follows from (b).

In any case, the corresponding topological group satisfies the assumptions of
Theorem \ref{mooregeneral} and we obtain the conclusion.
\end{pf}

{\today}

\end{document}